\documentstyle[graphicx]{article}[15pt]
\input amssym.def

\catcode`@=12
 \long\def\@makefntext#1{\noindent #1}
\newskip\tabcentering \tabcentering=1000pt plus 1000pt minus 1000pt

\def\MCH#1#2{\setbox0=\hbox{\raise#1\hbox{#2}}\smash{\box0}}

\def\@evenfoot{}\def\@oddfoot{}

\def\sec#1{\vspace{5mm}\leftline{\bf #1}\vspace{3mm}}

\catcode`@=12

\def\bc{\begin{center}}
\def\ec{\end{center}}

\def\hang{\hangindent\parindent}
\def\textindent#1{\indent\llap{\qquad #1\ \ \enspace}\ignorespaces}
\def\ref{\par\hang\textindent}

\def\a1{(a_1, a_2, \cdots, a_n)}

\def\a{\alpha}

\begin{document}
\thispagestyle{empty}
\vspace*{-3.0truecm}
\noindent
\vspace{1 true cm}
 \bc{\large\bf The   pressureless limits of  Riemann solutions to the   Euler equations of one-dimensional compressible
fluid flow
with a source term

\footnotetext{$^{*}$Corresponding author. \\
\indent \,\,\,\,\,\,\,\,E-mail address:  zqshao@fzu.edu.cn.\\
 }}\ec

 \vspace*{0.2 true cm}
\bc{\bf Shouqiong Sheng$^{a}$, Zhiqiang  Shao$^{a, *}$   \\
{\it $^{a}$College of Mathematics and Computer Science, Fuzhou University, Fuzhou 350108, China}
 }\ec

 \vspace*{2.5 true mm}
\setlength{\unitlength}{1cm}
\begin{picture}(20,0.1)
\put(-0.6,0){\line(1,0){14.5}}
\end{picture}

 \vspace*{2.5 true mm}
\noindent{\small {\small\bf Abstract}

 \vspace*{2.5 true mm}In this paper, we study  the  limits of
  Riemann solutions to the   inhomogeneous Euler equations of one-dimensional
 compressible fluid flow as the  adiabatic exponent $\gamma$ tends to one.
Different from the homogeneous equations, the Riemann solutions of the inhomogeneous system are non self-similar.  It is rigorously shown that,  as $\gamma$ tends to one, any two-shock Riemann solution
 tends to a delta shock solution  of the pressureless Euler system with a source term, and
the intermediate density between the two shocks tends to a weighted $\delta$-mesaure which forms the delta
shock; while any two-rarefaction-wave Riemann solution tends to a two-contact-discontinuity
solution  of the pressureless Euler system with a source term, whose intermediate state between the two contact discontinuities
is a vacuum state. Moreover, we also give some numerical results to confirm the theoretical analysis.

 \vspace*{2.5 true mm}
\noindent{\small {\small\bf MSC: } 35L65;  35L67

 \vspace*{2.5 true mm}
\noindent{\small {\small\bf Keywords:}     Pressureless limit;  Inhomogeneous Euler equations of one-dimensional compressible
fluid flow;   Non self-similar Riemann solution

 \vspace*{2.5 true mm}
\setlength{\unitlength}{1cm}
\begin{picture}(20,0.1)
\put(-0.6,0){\line(1,0){14.5}}
\end{picture}



\baselineskip 15pt
 \sec{\Large\bf 1.\quad  Introduction }
  The Euler equations of one-dimensional compressible fluid flow  with the Coulomb-like friction term can be written as
$$ \left\{\begin{array}{ll} \rho_{t}+(\rho u)_{x}=0,\\u_{t}+(\frac{u^{2}}{2}+p (\rho))_{x}=\beta,\end{array}\right .\eqno{(1.1)}
$$
where $\beta$ is a constant, the nonlinear   function  $p(\rho)=\frac{\theta}{2}\rho^{\gamma-1},$   $\theta=\frac{\gamma-1}{2}$  and  $\gamma\in(1,2)$   is a constant.

Shen [24] considered
 the pressureless Euler system with the Coulomb-like friction term and obtained the non self-similar Riemann solutions by introducing a new velocity:
$$v(t,x)=u(t,x)-\beta t,\eqno{(1.2)}$$
which was introduced by  Faccanoni and Mangeney [9] to study the Riemann problem of the shallow water equations with the Coulomb-like friction term.

If $\beta=0$, then
 the system (1.1) becomes the homogeneous  Euler equations of one-dimensional compressible fluid flow (cf. [8]):
$$ \left\{\begin{array}{ll} \rho_{t}+(\rho u)_{x}=0,\\u_{t}+(\frac{u^{2}}{2}+p (\rho))_{x}=0.\end{array}\right .\eqno{(1.3)}
$$

System (1.3) was firstly derived by Earnshaw [8] in 1858 for isentropic flow and is also viewed as
 the Euler equations of one-dimensional compressible fluid flow [14]. where $\rho$ denotes the density, $u$ the velocity, and $p(\rho)$ the
pressure of the fluid.
System (1.3) has other different physical backgrounds. For instance, it is a scaling
limit system of Newtonian dynamics with long-range interaction for a continuous
distribution of mass   in $R$ [20, 21] and also a hydrodynamic
limit for the Vlasov equation [1].

The  solutions for system (1.3)  were widely studied by many scholars (see [4-5, 7-8,  17-18, 22] ). In particular,  the existence of global weak solutions of the Cauchy problem  was first established by DiPerna [7] for the case of $1 <\gamma<3$ by using the Glimm's
scheme method.
 Using the result of DiPerna [7], Li [17] obtained a global weak solution to the Cauchy problem for
the case $  -1 < \gamma < 1.$  Using the theory of compensated compactness coupled with some
basic ideas of the kinetic formulation, Lu [18]  established an existence theorem for global entropy solutions for the case $\gamma > 3$. Cheng [5] also used the same methods as in [18]
 to obtain  the existence of global entropy solutions for the Cauchy
problem with a uniform amplitude bound for
the case $1 <\gamma  < 3.$

  When $\gamma\rightarrow 1$, the limiting system of (1.1) formally becomes  the pressureless Euler system with the Coulomb-like friction
term,
  $$ \left\{\begin{array}{ll} \rho_{t}+(\rho u)_{x}=0,\\u_{t}+(\frac{u^{2}}{2})_{x}=\beta,\end{array}\right .\eqno{(1.4)}
$$
which can  be  also obtained   by taking the constant pressure where the force is assumed to be the gravity with $\beta$  being the gravity constant [6].

 For the Euler system  of power law in Eulerian coordinates,
$$ \left\{\begin{array}{ll} \rho_{t}+(\rho u)_{x}=0,\\(\rho u)_{t}+
 (\rho u^{2}+p(\rho) )_{x}=0,\end{array}\right .
\eqno{(1.5)}$$  when the pressure tends to zero or a constant, the Euler system  (1.5)  formally tends to the zero pressure gas dynamics.  In earlier seminal papers,
Chen and Liu   [2] first showed the formation of $\delta$-shocks and vacuum states of the Riemann solutions to the Euler system (1.5) for polytropic gas by taking limit $\varepsilon \rightarrow 0+$ in the model $p(\rho)=\varepsilon\rho^{\gamma}/\gamma$     $ ( \gamma >1)$,
which describe the phenomenon of concentration and cavitation rigorously in mathematics.  Further, they also obtained the same results for the Euler
equations for nonisentropic
fuids in [3]. The same problem for the  Euler equations  (1.5) for
isothermal case $( \gamma=1)$ was studied by Li [16].  Recently,  Muhammad Ibrahim, Fujun Liu and Song Liu
[12] showed the same phenomenon of concentration also exists in the mode $p(\rho) = \rho^{\gamma}$ $(0 < \gamma < 1)$ as $\gamma \rightarrow 0$, which is the case that the pressure goes to a constant. Namely, they  showed rigorously the formation of delta wave with the limiting behavior of Riemann solutions to the Euler equations (1.5). For  some other physical models,  there are also many results, the readers are referred to [10, 11, 19,23, 25-27,  30-32] and
the references cited therein.

    Motivated by [2-3, 16], in this paper, we focus on the pressureless limits  of Riemann solutions to the inhomogeneous Euler system (1.1) of one-dimensional compressible fluid flow. Different from the homogeneous equations, the Riemann solutions  are non self-similar,  we show the same phenomenon of concentration and cavitation also exists in the case $1 < \gamma < 2$ as $\gamma \rightarrow1$.

\vskip 0.1in
The organization of this article is as follows: In section 2 and section 3, we display some results on  the Riemann solutions of (1.4), (1.1), respectively. In section 4, we show rigorously  the formation of $\delta$-shocks and  vacuum states in the   pressureless limit of  Riemann solutions to (1.1) as $\gamma \rightarrow1$.
In Section 5, we present some representative numerical results to demonstrate the validity of the theoretical analysis in Sections 4.

\baselineskip 15pt
 \sec{\Large\bf 2.\quad   Preliminaries }
  In this section, we give  the results on the Riemann problem for system (1.4).   For the homogeneous pressureless Euler system  corresponding to system (1.4),  the results on the Riemann problem
 can be found in [28, 26, 30, 13].

By a change of variable (1.2), system (1.4) can be  rewritten in the conservative form
$$ \left\{\begin{array}{ll} \rho_{t}+(\rho (v+\beta t))_{x}=0,\\v_{t}+(\frac{(v+\beta t)^{2}}{2})_{x}=0.\end{array}\right .\eqno{(2.1)}
$$

In this section, we are interested in
   the Riemann problem for (2.1)   with initial data
$$ (\rho, v)(0, x) =\left\{\begin{array}{ll} (\rho_{-},
u_{-}),\,\,\,\,x< 0,\\(\rho_{+},
u_{+}),\,\,\,\,x> 0,\end{array} \right.\eqno{(2.2)}$$
where  $\rho_{\pm}>0$ and $u_{\pm}$   are given constant states.

It can be seen that the  solutions of the Riemann problem to system (1.4) can be obtained from the corresponding ones of (2.1) and (2.2) by using the change of state  variables $(\rho, u)(t, x)
=(\rho, v+ \beta t)(t, x)$  directly.

The system (2.1)  has a double eigenvalue  $\lambda =v+\beta t$ whose corresponding right
eigenvector is $\overrightarrow{r} = (1, 0)^{T}.$ Since $\nabla\lambda\cdot \overrightarrow{r}\equiv 0,$ so (2.1) is full linear degenerate and elementary waves are contact discontinuities.

For a discontinuity $\sigma(t)=x'(t),$ the Rankine-Hugoniot conditions $$\left\{
    \begin{array}{ll}
      -\sigma(t)[\rho]+[\rho (v+\beta t)]=0, \\
     -\sigma(t)[v]+[\frac{(v+\beta t)^{2}}{2}]=0,
    \end{array}
  \right.\eqno (2.3)
$$
hold, where $[\rho]=\rho -\rho_{-},$   etc. By solving  (2.3), we obtain  contact discontinuity $J(\rho_{-},u_{-} ):$
$$\sigma(t)=v+\beta t=u_{-}+\beta t.\eqno (2.4)
$$

We now  can construct  the Riemann solutions of (2.1) and (2.2)  by contact discontinuities, vacuum  or $\delta$-shock wave connecting two constant states $(\rho_{\pm},u_{\pm})$.

For the case $u_{-}< u_{+}$, the Riemann solution   consists of two contact discontinuities with a vacuum between them,  which is shown as
$$ (\rho,v)(t,x)=\left\{\begin{array}{ll} (\rho_{-}, u_{-}),\,\,\,\,\,\,\,\,\,\,\,\,\,\,\,\,\,\,\,\,\,\,\,\,\,\,-\infty<x<u_{-}t+\frac{1}{2}\beta t^{2}, \\Vac,\,\,\,\,\,\,\,\,\,\,\,\,\,\,\,\,\,\,\,\,\,\,\,\,\,\,\,\,\,\,\,\,\,\,\,\,\,u_{-}t+\frac{1}{2}\beta t^{2} \leq x \leq u_{+}t+\frac{1}{2}\beta t^{2},\\(\rho_{+},u_{+}), \,\,\,\,\,\, \,\,\,\,\,\,\,\,\,\,\,\,\,\,\,\,\,\,\,\, u_{+}t+\frac{1}{2}\beta t^{2}<x<+\infty.
\end{array}\right .\eqno{(2.5)}
$$
The Riemann solution can be expressed by:
$$ (\rho_{-}, u_{-})+J_{1}+Vac+J_{2}+(\rho_{+},u_{+}),\eqno{(2.6)}
$$where   ``+" means ``followed by".

For the case $u_{-}= u_{+}$, the Riemann solution   consists of one contact discontinuity,  which is shown as
$$ (\rho,v)(t,x)=\left\{\begin{array}{ll} (\rho_{-}, u_{-}),\,\,\,\,\,\,\,\,\,\,\,\,\,\,\,\,\,\,\,\,\,\,\,\,\,\,-\infty<x<u_{-}t+\frac{1}{2}\beta t^{2}, \\(\rho_{+},u_{+}), \,\,\,\,\,\, \,\,\,\,\,\,\,\,\,\,\,\,\,\,\,\,\,\,\,\, u_{-}t+\frac{1}{2}\beta t^{2}<x<+\infty.
\end{array}\right .\eqno{(2.7)}
$$
The Riemann solution can be expressed by:
$$(\rho_{-}, u_{-})+ J + (\rho_{+},u_{+}).\eqno{(2.8)}
$$

For the case $u_{-}> u_{+}$, the Riemann solution cannot be constructed by using the classical waves,
and the delta shock wave appears. The Riemann solution can be expressed by:
$$(\rho_{-}, u_{-})+ \delta S + (\rho_{+},u_{+}).\eqno{(2.9)}
$$
The delta shock $\delta S$ satisfies the generalized Rankine-Hugoniot conditions
$$
\left\{
     \begin{array}{ll}
       \frac{dx(t)}{dt}=u_{\delta}(t), \\
       \frac{dw(t)}{dt}=u_{\delta}(t) [\rho]-[\rho(v+\beta t)], \\
       u_{\delta}(t) [v]=[\frac{(v+\beta t)^{2}}{2} ],
     \end{array}
   \right.\eqno{(2.10)}$$
where $[\rho] = \rho_{+}- \rho_{-}$, $x(t)$, $w(t)$ and $u_{\delta}(t) =v_{\delta}+\beta t$   respectively denote the location,  weight and propagation speed of the delta shock, and $(x, w)(0) =(0, 0).$

By simple calculation,
we have$$v_{\delta}=\frac{1}{2}(u_{-}+u_{+}), ~~~~x(t) =v_{\delta}t+\frac{1}{2}\beta t^{2},~~~~
 w(t)=\frac{1}{2}(\rho_{-}+\rho_{+})(u_{-}-u_{+})\,t. \eqno{(2.11)}$$

We also can justify that  the delta shock satisfies the generalized entropy  condition
$$u_{+}+\beta t < u_{\delta}(t) < u_{-}+\beta t. \eqno{(2.12)}$$
Thus, we have obtained the Riemann solutions of (2.1) and (2.2).

In summary, we obtain the Riemann solutions to system (1.4) as follows

(1) For $u_{-} >u_{+}$, the Riemann solution to system (1.4) has
the following form:
$$(\rho, u)(t,x)=\left\{
                  \begin{array}{ll}
                    (\rho_{-}, u_{-}+\beta t),& \hbox{$x<x(t)$,}\\(w(t)\delta(x-x(t)),u_{\delta} (t)), & \hbox{$x=x(t)$,}\\
                    (\rho_{+}, u_{+}+\beta t), & \hbox{$x>x(t)$,}
                  \end{array}
                \right.
\eqno{(2.13)}$$
where
$$x(t)=\frac{1}{2}(u_{-}+u_{+})t+\frac{1}{2}\beta t^{2},~~~~
 w(t)=\frac{1}{2}(\rho_{-}+\rho_{+})(u_{-}-u_{+})\,t, ~~~~u_{\delta}(t)=\frac{1}{2}(u_{-}+u_{+})+\beta t.\eqno{(2.14)}$$

(2) For $u_{-}<u_{+}$, the Riemann solution can be expressed as
 $$ (\rho,u)(t,x)=\left\{\begin{array}{ll} (\rho_{-}, u_{-}+\beta t),\,\,\,\,\,\,\,\,\,\,\,\,\,\,\,\,\,\,\,\,\,\,\,\,\,\,-\infty<x<u_{-}t+\frac{1}{2}\beta t^{2}, \\Vac,\,\,\,\,\,\,\,\,\,\,\,\,\,\,\,\,\,\,\,\,\,\,\,\,\,\,\,\,\,\,\,\,\,\,\,\,\,u_{-}t+\frac{1}{2}\beta t^{2} \leq x \leq u_{+}t+\frac{1}{2}\beta t^{2},\\(\rho_{+},u_{+}+\beta t), \,\,\,\,\,\, \,\,\,\,\,\,\,\,\,\,\,\,\,\,\,\,\,\,\,\, u_{+}t+\frac{1}{2}\beta t^{2}<x<+\infty,
\end{array}\right .\eqno{(2.15)}
$$
 where  the locations and propagation speeds of two contact discontinuities $J_{1}$ and $J_{2}$ are identical with those in the Riemann
solution of (2.1) and (2.2).

(3) For  $u_{-}=u_{+}$, the Riemann solution can be expressed as
 $$ (\rho,u)(t,x)=\left\{\begin{array}{ll} (\rho_{-}, u_{-}+\beta t),\,\,\,\,\,\,\,\,\,\,\,\,\,\,\,\,\,\,\,\,\,\,\,\,\,\,-\infty<x<u_{-}t+\frac{1}{2}\beta t^{2}, \\(\rho_{+},u_{+}+\beta t), \,\,\,\,\,\, \,\,\,\,\,\,\,\,\,\,\,\,\,\,\,\,\,\,\,\, u_{-}t+\frac{1}{2}\beta t^{2}<x<+\infty,
\end{array}\right .\eqno{(2.16)}
$$
 where  the location and propagation speed of  contact discontinuity $J_{}$  are identical with those in the Riemann
solution of (2.1) and (2.2).

\baselineskip 15pt
 \sec{\Large\bf 3.\quad   Riemann problem   for  Euler equations  with a source term (1.1)}

 \indent

 In this section, we  construct the
Riemann solutions of the Euler equations  with the Coulomb-like friction term
(1.1).

Using (1.2), system (1.1) is rewritten in the conservative form
$$ \left\{\begin{array}{ll} \rho_{t}+(\rho (v+\beta t))_{x}=0,\\v_{t}+(\frac{(v+\beta t)^{2}}{2}+\frac{\gamma-1}{4}\rho^{\gamma-1} )_{x}=0.\end{array}\right .\eqno{(3.1)}
$$

In this section, we are interested in
   the Riemann problem for (3.1)   with initial data
$$ (\rho, v)(0, x) =\left\{\begin{array}{ll} (\rho_{-},
u_{-}),\,\,\,\,x< 0,\\(\rho_{+},
u_{+}),\,\,\,\,x> 0,\end{array} \right.\eqno{(3.2)}$$
where  $\rho_{\pm}>0$ and $u_{\pm}$   are given constant states.

The system (3.1) can be  reformulated  in a quasi-linear form
$$\left(\begin{array}{cc}\rho\\v
 \end{array}\right)_t+\left(
  \begin{array}{ccc}v+\beta t& \rho \\\frac{(\gamma-1)^{2}}{4}\rho^{\gamma-2} & v+\beta t
  \end{array}\right)\left(\begin{array}{cc}\rho\\v
 \end{array}\right)_x=\left(\begin{array}{cc}0\\0
 \end{array}\right).\eqno{(3.3)}
$$ By (3.3), it is easy  to see that    system (3.1) has two  eigenvalues
 $$ \lambda_{1}^{\gamma}=v+\beta t-\frac{\gamma-1}{2}\rho^{\frac{\gamma-1}{2}},\,\,\,\,\,\,\,\lambda_{2}^{\gamma}=v+\beta t+\frac{\gamma-1}{2}\rho^{\frac{\gamma-1}{2}},\eqno{(3.4)}
$$
with the corresponding right eigenvectors
$$\overrightarrow{r}_{1}^{\gamma} =(1, -\frac{\gamma-1}{2}\rho^{\frac{\gamma-3}{2}})^{T}, \,\,
\overrightarrow{r}_{2}^{\gamma} =(1, \frac{\gamma-1}{2}\rho^{\frac{\gamma-3}{2}})^{T},
 $$satisfying
$$\nabla\lambda_{1}^{\gamma}\cdot \overrightarrow{r}_{1}^{\gamma}=-\frac{(\gamma-1)(\gamma+1)}{4}\rho^{\frac{\gamma-3}{2}}<0,$$ $$\nabla\lambda_{2}^{\gamma}\cdot \overrightarrow{r}_{2}^{\gamma}=\frac{(\gamma-1)(\gamma+1)}{4}\rho^{\frac{\gamma-3}{2}}>0.$$Therefore,  system (3.1)  is strictly hyperbolic for $\rho>0$,
 both   characteristic fields  are genuinely nonlinear  and the associated waves are  shock waves  or rarefaction waves.

The Riemann invariants may be selected as
$$w^{\gamma} =v+\rho^{\frac{\gamma-1}{2}}, \,\,\,\,z^{\gamma} =v-\rho^{\frac{\gamma-1}{2}},\eqno{(3.5)}$$
which    satisfy $\bigtriangledown w^{\gamma} \cdot \overrightarrow{r_1}^{\gamma}=0$  and $\bigtriangledown z^{\gamma} \cdot \overrightarrow{r_2}^{\gamma}=0$, respectively.

Given a state $(\rho_{-}, u_{-})$, the rarefaction wave curves in the phase plane,
which are the sets of states that can be connected on the right by a 1-rarefaction
or 2-rarefaction wave, are as follows
$$R_{1}^{\gamma}(\rho_{-},u_{-}):\,\,\left\{\begin{array}{ll} \frac{dx}{dt}=\lambda_{1}^{\gamma}=v+\beta t-\frac{\gamma-1}{2}\rho^{\frac{\gamma-1}{2}},
\\v+\rho^{\frac{\gamma-1}{2}}=u_{-}+\rho_{-}^{\frac{\gamma-1}{2}},
\,\,\,\,\rho<\rho_{-}, v>u_{-},\\\lambda_{1}^{\gamma}(\rho_{-},u_{-})<\lambda_{1}^{\gamma}(\rho,v),\end{array} \right.\eqno{(3.6)}$$ and
$$R_{2}^{\gamma}(\rho_{-},u_{-}):\,\,\left\{\begin{array}{ll}  \frac{dx}{dt}=\lambda_{2}^{\gamma}=v+\beta t+\frac{\gamma-1}{2}\rho^{\frac{\gamma-1}{2}},
\\v-\rho^{\frac{\gamma-1}{2}}=u_{-}-\rho_{-}^{\frac{\gamma-1}{2}},
\,\,\,\,\rho>\rho_{-}, v>u_{-},\\\lambda_{2}^{\gamma}(\rho_{-},u_{-})<\lambda_{2}^{\gamma}(\rho,v).\end{array} \right.\eqno{(3.7)}$$

Differentiating $v$ with respect to $\rho$
in the second equation of  (3.6), we have $$\frac{dv}{d\rho} = -\frac{\gamma-1}{2}\rho^{\frac{\gamma-3}{2}}<0,$$
$$\frac{d^{2}v}{d\rho^{2}} =-\frac{(\gamma-1)(\gamma-3)}{4}\rho^{\frac{\gamma-5}{2}}>0,$$
which implies that the 1-rarefaction wave curve $R_{1}^{\gamma}(\rho_{-},u_{-})$ is monotonic decreasing and convex in the $(\rho, v)$ phase plane. Similarly, one can also obtain $\frac{dv}{d\rho}> 0$ and $\frac{d^{2}v}{d\rho^{2}} < 0$ by differentiating $v$ with respect to $\rho$
in the second equation of  (3.7),  which implies that
the 2-rarefaction wave curve $R_{2}^{\gamma}(\rho_{-},u_{-})$ is monotonic increasing and concave in the $(\rho, v)$ phase plane.  Moreover, it  can be concluded  from (3.6) that
 $\lim\limits_{\rho\rightarrow 0^{+}}v= u_{-}+\rho_{-}^{\frac{\gamma-1}{2}}$  for the 1-rarefaction wave curve $R_{1}^{\gamma}(\rho_{-},u_{-})$, which indicates that curve
$R_{1}^{\gamma}(\rho_{-},u_{-})$ intersects the $v$-axis at the point $(0,\widetilde{v}_{\ast}^{\gamma})$,  where $\widetilde{v}_{\ast}^{\gamma}$ is determined by
$\widetilde{v}_{\ast}^{\gamma}= u_{-}+\rho_{-}^{\frac{\gamma-1}{2}}$. It can also be  seen from (3.7) that $\lim\limits_{\rho\rightarrow +\infty}v=+\infty$
 for the 2-rarefaction
wave curve $R_{2}^{\gamma}(\rho_{-},u_{-})$.

  Let $\sigma^{\gamma}(t)=\frac{dx^{\gamma} (t)}{dt}$ be the speed of a bounded discontinuity $x=x^{\gamma}(t)$, then the Rankine-Hugoniot conditions for the conservative system (3.1) are given by
  $$\left\{
    \begin{array}{ll}
      -\sigma^{\gamma}(t)[\rho]+[\rho (v+\beta t)]=0, \\
     -\sigma^{\gamma}(t)[ v]+[\frac{(v+\beta t)^{2}}{2}+\frac{\gamma-1}{4}\rho^{\gamma-1}]=0,
    \end{array}
  \right.\eqno (3.8)
$$
 where $[\rho]=\rho-\rho_{-}$, etc. From (3.8) we have
$$
      \sigma^{\gamma} (t)=\frac{[\rho  (v+\beta t)]}{[\rho]},$$
     $$\frac{v-v_{-}}{\rho-\rho_{-}} =\pm\sqrt{\frac{\frac{\gamma-1}{2}[\rho^{\gamma-1}]}{(\rho+\rho_{-})[\rho]}},
     \eqno (3.9)$$
where  $(\rho_{-},v_{-})$ and $(\rho_{},v_{})$ are  the left state and the right state, respectively.

1-$shock$  $curve$  $S_{1}^{\gamma} (\rho_{-}, u_{-})$:

The  Lax
 entropy condition  implies  that the propagation speed $\sigma_{1}^{\gamma} (t)$ for the 1-shock wave $S_{1}^{\gamma} $
has to be satisfied with$$\sigma_{1}^{\gamma} (t)<\lambda_{1}^{\gamma} (\rho_{-},v_{-}),\,\,\,\,\lambda_{1}^{\gamma} (\rho,v)<\sigma_{1}^{\gamma} (t)<\lambda_{2}^{\gamma} (\rho,v).\eqno (3.10)
$$
From the first equation of (3.8), we obtain
$$\sigma_{1}^{\gamma} (t)=\frac{\rho (v+\beta t) -\rho_{-}(v_{-}+\beta t)}{\rho-\rho_{-}} = v_{-}+\beta t+\frac{\rho}{\rho-\rho_{-}}(v-v_{-}).\eqno (3.11)$$
Then,    substituting  (3.11) into the  first inequality  of (3.10), we have
$$\frac{\rho}{\rho-\rho_{-}}(v-v_{-})<-\frac{\gamma-1}{2}\rho_{-}^{\frac{\gamma-1}{2}}<0,$$
which shows that $v - v_{-}$ and $\rho- \rho_{-}$ have different signs. Thus,
from (3.9)  we have $$v=v_{-}-\sqrt{\frac{\frac{\gamma-1}{2}(\rho^{\gamma-1}-\rho_{-}^{\gamma-1})}{(\rho+\rho_{-})(\rho-\rho_{-})}}(\rho-\rho_{-}).
$$
 If $v>v_{-}$, then $\rho<\rho_{-}$, and
 $$\sigma_{1}^{\gamma} (t)- v_{-}-\beta t=\frac{\rho}{\rho-\rho_{-}}(v-v_{-})=-\rho\sqrt{\frac{\frac{\gamma-1}{2}(\rho^{\gamma-1}-\rho_{-}^{\gamma-1})}{(\rho+\rho_{-})(\rho-\rho_{-})}}=-\frac{\gamma-1}{2}
 \overline{\rho}^{\frac{\gamma-2}{2}}\rho\sqrt{\frac{2}{\rho+\rho_{-}}},
$$
for some $\bar{\rho}\in(\rho,\rho_{-}).$  By direct calculation, we have
$$\frac{\gamma-1}{2}\rho_{-}^{\frac{\gamma-1}{2}}-\frac{\gamma-1}{2}
 \overline{\rho}^{\frac{\gamma-2}{2}}\rho\sqrt{\frac{2}{\rho+\rho_{-}}}>\frac{\gamma-1}{2}\bigg(\rho_{-}^{\frac{\gamma-1}{2}}-
 \rho^{\frac{\gamma-2}{2}}\rho\sqrt{\frac{2}{\rho+\rho_{-}}}\bigg)>\frac{\gamma-1}{2}(\rho_{-}^{\frac{\gamma-1}{2}}-\rho^{\frac{\gamma-1}{2}})>0,
$$ which implies that
$$\sigma_{1}^{\gamma} (t)- v_{-}-\beta t >-\frac{\gamma-1}{2}\rho_{-}^{\frac{\gamma-1}{2}}.
$$
This  contradicts with $\sigma_{1}^{\gamma} (t)<\lambda_{1}^{\gamma} (\rho_{-},v_{-})$.  Hence, given a state $(\rho_{-},u_{-})$,  the 1-shock wave curve $S_{1}^{\gamma} (\rho_{-},u_{-})$
  in the phase plane which is the set of states that can be connected on the right by a 1-shock is  as follows
$$S_{1}^{\gamma} (\rho_{-},u_{-}):\,\,\left\{\begin{array}{ll} \sigma_{1}^{\gamma} (t)=u_{-}
+\beta t-\rho\sqrt{\frac{\frac{\gamma-1}{2}(\rho^{\gamma-1}-\rho_{-}^{\gamma-1})}{(\rho+\rho_{-})(\rho-\rho_{-})}},
\\v=u_{-}-\sqrt{\frac{\frac{\gamma-1}{2}(\rho^{\gamma-1}-\rho_{-}^{\gamma-1})}{(\rho+\rho_{-})(\rho-\rho_{-})}}(\rho-\rho_{-}),~~~\\\rho>\rho_{-},v<u_{-}.\end{array} \right.\eqno{(3.12)}$$

2-$shock~curve~S_{2}^{\gamma} (\rho_{-},u_{-})$:

Similarly, the propagation speed $\sigma_{2}^{\gamma} (t)$ for the 2-shock wave $S_{2}^{\gamma} $ should satisfy
$$\lambda_{1}^{\gamma} (\rho_{-},v_{-})<\sigma_{2}^{\gamma} (t)<\lambda_{2}^{\gamma} (\rho_{-},v_{-}),\,\,\,\,\lambda_{2}^{\gamma} (\rho,v)<\sigma_{2}^{\gamma} (t).
$$Then,  given a state $(\rho_{-},u_{-})$,  the 2-shock wave curve $S_{2}^{\gamma} (\rho_{-},u_{-})$
  in the phase plane which is the set of states that can be connected on the right by a 2-shock is  as follows
$$S_{2}^{\gamma} (\rho_{-},u_{-}):\,\,\left\{\begin{array}{ll} \sigma_{2}^{\gamma} (t)=u_{-}
+\beta t+\rho\sqrt{\frac{\frac{\gamma-1}{2}(\rho^{\gamma-1}-\rho_{-}^{\gamma-1})}{(\rho+\rho_{-})(\rho-\rho_{-})}},
\\v=u_{-}
+\sqrt{\frac{\frac{\gamma-1}{2}(\rho^{\gamma-1}-\rho_{-}^{\gamma-1})}{(\rho+\rho_{-})(\rho-\rho_{-})}}(\rho-\rho_{-}),~~~\\\rho<\rho_{-},v<u_{-}.\end{array} \right.\eqno{(3.13)}$$

Differentiating $v$ with respect to $\rho$
 in the second equation in  (3.12)  yields that for $\rho>\rho_{-}$,
$$\frac{dv}{d\rho}=-\frac{1}{2}\sqrt{\frac{\frac{\gamma-1}{2}(\rho+\rho_{-})}{(\rho^{\gamma-1}-\rho_{-}^{\gamma-1})(\rho-\rho_{-})}}
\frac{(\gamma-1)\rho^{\gamma-2}(\rho-\rho_{-})(\rho+\rho_{-})+2\rho_{-}(\rho^{\gamma-1}-\rho_{-}^{\gamma-1})}{(\rho+\rho_{-})^{2}}<0,$$
which indicates that the 1-shock wave curve $S_{1}^{\gamma} (\rho_{-},u_{-})$ is monotonic decreasing in the $(\rho, v) $ phase plane  ($\rho>\rho_{-})$. Similarly,
from  (3.13), for $\rho<\rho_{-}$ we have $\frac{dv}{d\rho} > 0, $  which indicates that the 2-shock wave
curve $S_{2}^{\gamma} (\rho_{-},u_{-})$ is monotonic increasing in the $(\rho, v) $ phase plane  ($\rho<\rho_{-})$.  It can be seen from (3.13)
that $\lim\limits_{\rho\rightarrow 0^{+}}v= u_{-}-\sqrt{\frac{\gamma-1}{2}}\rho_{-}^{\frac{\gamma-1}{2}}$ for the 2-shock wave curve $S_{2}^{\gamma} (\rho_{-},u_{-})$, which implies that $S_{2}^{\gamma} (\rho_{-},u_{-})$ intersects the $v$-axis at the point $(0,\widetilde{v}_{\ast\ast}^{\gamma}   )$,  where $\widetilde{v}_{\ast\ast}^{\gamma} $
is determined by  $\widetilde{v}_{\ast\ast}^{\gamma} =u_{-}-\sqrt{\frac{\gamma-1}{2}}\rho_{-}^{\frac{\gamma-1}{2}}.$
 It can also be derived from (3.12) that $\lim\limits_{\rho\rightarrow +\infty}v= -\infty$
 for the 1-shock wave curve $S_{1}^{\gamma} (\rho_{-},u_{-})$.

In the $(\rho, v) $ phase plane, through a given point $(\rho_{-}, u_{-})$, we draw the elementary wave curves  $R_{j}^{\gamma} (\rho_{-}, u_{-})$  and $S_{j}^{\gamma} (\rho_{-}, u_{-})$ (j=1, 2).  These  elementary wave curves divide the  $(\rho, v)$ phase plane into five regions (see Fig. 1).  According to the right
state $(\rho_{+},u_{+})$ in the different regions, one can construct the unique global Riemann solution of (3.1)  and (3.2) as follows:

(1) $(\rho_{+},u_{+})\in I(\rho_{-},u_{-}):$ $(\rho_{-},u_{-})+R_{1}^{\gamma} +(\rho_{\ast\gamma},v_{\ast\gamma})+R_{2}^{\gamma} +(\rho_{+},u_{+});$

 (2)$(\rho_{+},u_{+})\in II(\rho_{-},u_{-}):$ $(\rho_{-},u_{-})+S_{1}^{\gamma} +(\rho_{\ast\gamma},v_{\ast\gamma})+R_{2}^{\gamma} +(\rho_{+},u_{+});$

 (3)$(\rho_{+},u_{+})\in III(\rho_{-},u_{-}):$ $(\rho_{-},u_{-})+R_{1}^{\gamma} +(\rho_{\ast\gamma},v_{\ast\gamma})+S_{2}^{\gamma} +(\rho_{+},u_{+});$

 (4)$(\rho_{+},u_{+})\in IV(\rho_{-},u_{-}):$ $(\rho_{-},u_{-})+S_{1}^{\gamma} +(\rho_{\ast\gamma},v_{\ast\gamma})+S_{2}^{\gamma} +(\rho_{+},u_{+});$

  (5)$(\rho_{+},u_{+})\in V(\rho_{-},u_{-}):$ $(\rho_{-},u_{-})+R_{1}^{\gamma} +\mathrm{Vac}+R_{2}^{\gamma} +(\rho_{+},u_{+}),$
\\ where $(\rho_{\ast\gamma},v_{\ast\gamma})$
is the intermediate state. By using (1.2), we obtain the Riemann solutions of (1.1) as follows

(1) $(\rho_{+},u_{+})\in I(\rho_{-},u_{-}):$ $(\rho_{-},u_{-}+\beta t)+R_{1}^{\gamma} +(\rho_{\ast\gamma},v_{\ast\gamma}+\beta t)+R_{2}^{\gamma} +(\rho_{+},u_{+}+\beta t);$

 (2)$(\rho_{+},u_{+})\in II(\rho_{-},u_{-}):$ $(\rho_{-},u_{-}+\beta t)+S_{1}^{\gamma} +(\rho_{\ast\gamma},v_{\ast\gamma}+\beta t)+R_{2}^{\gamma} +(\rho_{+},u_{+}+\beta t);$

 (3)$(\rho_{+},u_{+})\in III(\rho_{-},u_{-}):$ $(\rho_{-},u_{-}+\beta t)+R_{1}^{\gamma} +(\rho_{\ast\gamma},v_{\ast\gamma}+\beta t)+S_{2}^{\gamma} +(\rho_{+},u_{+}+\beta t);$

 (4)$(\rho_{+},u_{+})\in IV(\rho_{-},u_{-}):$ $(\rho_{-},u_{-}+\beta t)+S_{1}^{\gamma} +(\rho_{\ast\gamma},v_{\ast\gamma}+\beta t)+S_{2}^{\gamma} +(\rho_{+},u_{+}+\beta t);$

  (5)$(\rho_{+},u_{+})\in V(\rho_{-},u_{-}):$ $(\rho_{-},u_{-}+\beta t)+R_{1}^{\gamma} +\mathrm{Vac}+R_{2}^{\gamma} +(\rho_{+},u_{+}+\beta t).$

\hspace{65mm}\setlength{\unitlength}{0.8mm}\begin{picture}(80,66)
\put(-50,18){\vector(0,2){35}}
 \put(-48,0){\vector(2,0){103}}  \put(-53,49){$\rho$}
\put(56,-1){$v$}
\put(-38,7){$S_{2}^{\gamma}$}\put(6,5){}
\put(-32,48){$S_{1}^{\gamma}$}\put(32,48){$R_{2}^{\gamma}$}\put(42,19){$R_{2}^{\gamma}$}\put(12,5){$R_{1}^{\gamma}$}\put(45,6){V }\put(25,-5){$\widetilde{v}_{\ast}^{\gamma}$ }\put(-40,-5){$\widetilde{v}_{\ast\ast}^{\gamma}$ }
\put(1,12){$(\rho_{-},
u_{-})$}
\put(-5,6){III }\put(44,6){}\put(-5,29){II}
\put(-36,20){IV}\put(26,20){I}
\qbezier(-40,0)(14,12)(42,52)\qbezier(25,0)(-14,12)(-42,52)\qbezier(25,0)(36,06)(48,18)
\end{picture}
\vspace{0.6mm}  \vskip 0.2in \centerline{\bf Fig. 1.\,\,    Curves of  elementary waves.
   } \vskip 0.1in \indent

\baselineskip 15pt
 \sec{\Large\bf 4.\quad    Limits of Riemann solutions to (1.1)}In this section, we  study the limiting behavior   of the Riemann solutions to system (1.1)   as $\gamma$ tends to one, that is, the formation of delta shock and the vacuum states as $\gamma$ tends to one, respectively in the case $u_{+}<u_{-}$ and in the case $u_{+}>u_{-}$.

 \baselineskip 15pt
 \sec{\Large\bf 4.1.\quad    Formation of delta shock wave for  system (1.1) }In this subsection,  we study the phenomenon of the
concentration and the formation of delta shock in the Riemann solutions
to (1.1) in the case $u_{+}<u_{-}$  as $\gamma$ tends to one.

 \vskip 0.1in
\noindent{\small {\small\bf Lemma 4.1.} If $u_{+}<u_{-}$, then there is a sufficiently small $\gamma_{0} > 0$ such that
$(\rho_{+}, u_{+})\in IV( \rho_{-}, u_{-})$ as $1<\gamma <1 +\gamma_{0}$.

\vskip 0.1in
\noindent{\small {\small\bf Proof.} If $\rho_{+} = \rho_{-}$, then $(\rho_{+}, u_{+})\in IV( \rho_{-}, u_{-})$  for any $\gamma\in (1,2)$. Thus, we only need to consider the
case $\rho_{+} \neq\rho_{-}$.

By  (3.12) and (3.13), it is easy to see that all possible states $(\rho, v)$ that can be connected to the
left state $(\rho_{-}, u_{-})$ on the right by a 1-shock wave $S_{1}^{\gamma}$ or a 2-shock wave $S_{2}^{\gamma}$  satisfy
$$   S_{1}^{\gamma}: \,\,\,\,v=u_{-}-\sqrt{\frac{\frac{\gamma-1}{2}(\rho^{\gamma-1}-\rho_{-}^{\gamma-1})}{(\rho+\rho_{-})(\rho-\rho_{-})}}(\rho-\rho_{-}),~~~\rho>\rho_{-},\eqno{(4.1)}
$$
 $$S_{2}^{\gamma}: \,\,\,\,v=u_{-}+\sqrt{\frac{\frac{\gamma-1}{2}(\rho^{\gamma-1}-\rho_{-}^{\gamma-1})}{(\rho+\rho_{-})(\rho-\rho_{-})}}(\rho-\rho_{-}),~~  \rho<\rho_{-}.\eqno{(4.2)}
 $$
If $\rho_{+} \neq\rho_{-}$ and $(\rho_{+}, u_{+})\in IV( \rho_{-}, u_{-})$, then from Fig. 1, (4.1) and (4.2), we have$$  u_{+}<u_{-}-\sqrt{\frac{\frac{\gamma-1}{2}(\rho_{+}^{\gamma-1}-\rho_{-}^{\gamma-1})}{(\rho_{+}+\rho_{-})(\rho_{+}-\rho_{-})}}(\rho_{+}-\rho_{-}),~~~\mathrm {}  \,\,\rho_{+}>\rho_{-},\eqno{(4.3)}
$$
$$  u_{+}<u_{-}+\sqrt{\frac{\frac{\gamma-1}{2}(\rho_{+}^{\gamma-1}-\rho_{-}^{\gamma-1})}{(\rho_{+}+\rho_{-})(\rho_{+}-\rho_{-})}}(\rho_{+}-\rho_{-}),~~~\mathrm {}  \,\,\rho_{+}<\rho_{-}.\eqno{(4.4)}
$$
From (4.3) and (4.4), we derive that
$$ \sqrt{\frac{\frac{\gamma-1}{2}(\rho_{+}^{\gamma-1}-\rho_{-}^{\gamma-1})}{\rho_{+}^{2}-\rho_{-}^{2}}}< \frac{u_{-}-u_{+}}{|\rho_{+}-\rho_{-}|}. \eqno{(4.5)}$$
Since
$$  \lim_{{\gamma\rightarrow1}}\sqrt{\frac{\frac{\gamma-1}{2}(\rho_{+}^{\gamma-1}-\rho_{-}^{\gamma-1})}{\rho_{+}^{2}-\rho_{-}^{2}}}=0,
  \eqno{(4.6)}$$
it follows that there exists $\gamma_{0}>0$ small enough such that, when $ 1<\gamma <1 +\gamma_{0}$,  we have
$$ \sqrt{\frac{\frac{\gamma-1}{2}(\rho_{+}^{\gamma-1}-\rho_{-}^{\gamma-1})}{\rho_{+}^{2}-\rho_{-}^{2}}}< \frac{u_{-}-u_{+}}{|\rho_{+}-\rho_{-}|}. \eqno{}$$
  Then,  it is obvious
that $(\rho_{+}, u_{+})\in IV( \rho_{-}, u_{-})$ when $ 1<\gamma <1 +\gamma_{0}$.
 The proof is completed. $~\Box$

When  $1<\gamma <1+\gamma_{0}$, i.e., $(\rho_{+}, u_{+})\in IV( \rho_{-}, u_{-})$, suppose that $(\rho_{\ast\gamma},
v_{\ast\gamma})$ is the intermediate state
connected with $(\rho_{-}, u_{-})$ by a 1-shock wave $S_{1}^{\gamma}$ with the speed $\sigma_{1}^{\gamma}(t)$, and $(\rho_{+}, u_{+})$ by a 2-shock wave $S_{2}^{\gamma}$  with the speed $\sigma_{2}(t),$
 then it follows $$S_{1}^{\gamma} :\,\,\left\{\begin{array}{ll} \sigma_{1}^{\gamma} (t)=u_{-}
+\beta t-\rho_{\ast\gamma}\sqrt{\frac{\frac{\gamma-1}{2}(\rho_{\ast\gamma}^{\gamma-1}-\rho_{-}^{\gamma-1})}{(\rho_{\ast\gamma}+\rho_{-})(\rho_{\ast\gamma}-\rho_{-})}},
\\v_{\ast\gamma}=u_{-}-\sqrt{\frac{\frac{\gamma-1}{2}(\rho_{\ast\gamma}^{\gamma-1}-\rho_{-}^{\gamma-1})}{(\rho_{\ast\gamma}+\rho_{-})(\rho_{\ast\gamma}
-\rho_{-})}}(\rho_{\ast\gamma}-\rho_{-}),~~~\rho_{\ast\gamma}>\rho_{-},\end{array} \right.\eqno{(4.7)}$$$$S_{2}^{\gamma} :\,\,\left\{\begin{array}{ll} \sigma_{2}^{\gamma} (t)=v_{\ast\gamma}
+\beta t+\rho_{+}\sqrt{\frac{\frac{\gamma-1}{2}(\rho_{+}^{\gamma-1}-\rho_{\ast\gamma}^{\gamma-1})}{(\rho_{+}+\rho_{\ast\gamma})(\rho_{+}-\rho_{\ast\gamma})}},
\\u_{+}=v_{\ast\gamma}
+\sqrt{\frac{\frac{\gamma-1}{2}(\rho_{+}^{\gamma-1}-\rho_{\ast\gamma}^{\gamma-1})}{(\rho_{+}+\rho_{\ast\gamma})(\rho_{+}-\rho_{\ast\gamma})}}
(\rho_{+}-\rho_{\ast\gamma}),~~~\rho_{\ast\gamma}>\rho_{+}.\end{array} \right.\eqno{(4.8)}$$ From  (4.7) and (4.8), we have
$$ u_{-}-u_{+}=\sqrt{\frac{\frac{\gamma-1}{2}(\rho_{\ast\gamma}^{\gamma-1}-\rho_{-}^{\gamma-1})}{(\rho_{\ast\gamma}+\rho_{-})(\rho_{\ast\gamma}-\rho_{-})}}(\rho_{\ast\gamma}-\rho_{-})+
\sqrt{\frac{\frac{\gamma-1}{2}(\rho_{+}^{\gamma-1}-\rho_{\ast\gamma}^{\gamma-1})}{(\rho_{\ast\gamma}+\rho_{+})(\rho_{+}-\rho_{\ast\gamma})}}(\rho_{\ast\gamma}-\rho_{+}), \,\,\,\,   \rho_{\ast\gamma}>\rho_{\pm}.\eqno{(4.9)}
$$

 Then we have  the following lemmas.

 \vskip 0.1in
\noindent{\small {\small\bf Lemma 4.2.}
$\lim\limits_{\gamma\rightarrow1}\rho_{\ast\gamma}=+\infty,$ and $ \lim\limits_{\gamma\rightarrow1}\frac{\gamma-1}{2}\rho_{\ast\gamma}^{\gamma-1}=:a
=\frac{(u_{-}-u_{+})^{2}}{4}$.

\vskip 0.1in
\noindent{\small {\small\bf Proof.} Let $ \lim\limits_{\gamma\rightarrow1}\inf\rho_{\ast\gamma}=\alpha$, and $\lim\limits_{\gamma\rightarrow1}\sup\rho_{\ast\gamma}=\beta$.

If $ \alpha<\beta$ , then by the continuity of $\rho_{\ast\gamma}$, there exists a sequence $  \{\gamma_{n}\}_{n=1}^{\infty}\subseteq(1,2)$
such that
$$ \lim_{n\rightarrow +\infty}\gamma_{n}=1,\,\,\mathrm{ and}\,\, \lim_{n\rightarrow +\infty}\rho_{\ast\gamma_{n}}=c,$$
for  some $ c\in(\alpha,\beta).$ Then substituting the sequence  into the right-hand side of (4.9),  and  taking the limit $n\rightarrow +\infty$,
 we have
$$ \lim_{n\rightarrow+\infty}\sqrt{\frac{\frac{\gamma_{n}-1}{2}(\rho_{\ast\gamma_{n}}^{\gamma_{n}-1}-\rho_{\pm}^{\gamma_{n}-1})}
{\rho_{\ast\gamma_{n}}^{2}-\rho^{2}_{\pm}}}(\rho_{\ast\gamma_{n}}-\rho_{\pm})=0.
\eqno{(4.10)} $$
Thus,  we can obtain from (4.9) that
$$ u_{-}-u_{+}=0,$$
which contradicts with  $u_{-}>u_{+}$.
Then we must have $\alpha=\beta$, which   implies $\lim\limits_{\gamma\rightarrow1}\rho_{\ast\gamma}=\alpha.$

If  $\alpha\in(0,+\infty),$ then  we  can also get a contradiction when taking limit in (4.9). Thus $\alpha=0 $ or $ \alpha=+\infty$. By the condition
$\rho_{\ast\gamma}>\max\{\rho_{-},\rho_{+}\}$, it is easy to see that $\lim\limits_{\gamma\rightarrow1}\rho_{\ast\gamma}=\alpha=+\infty.$

Next taking the limit $ \gamma\rightarrow 1$ at the  right-hand side of (4.9), we have
$$\lim_{\gamma\rightarrow 1}\sqrt{\frac{\frac{\gamma-1}{2}(\rho_{\ast\gamma}^{\gamma-1}-\rho_{\pm}^{\gamma-1})}
{\rho_{\ast\gamma}^{2}-\rho_{\pm}^{2}}}(\rho_{\ast\gamma}-\rho_{\pm})
=\lim_{\gamma\rightarrow1}\sqrt{\frac{(\frac{\gamma-1}{2}\rho_{\ast\gamma}^{\gamma-1}-\frac{\gamma-1}{2}\rho_{\pm}^{\gamma-1})(\rho_{\ast\gamma}-\rho_{\pm})^{2}}
{\rho_{\ast\gamma}^{2}-\rho_{\pm}^{2}}}
=:\sqrt{a},$$
and
$$u_{-}-u_{+}=2\sqrt{a},\eqno{}
$$
from which we can get
$ a=\frac{(u_{-}-u_{+})^{2}}{4}.$
 The proof is completed. $~\Box$

\vskip 0.1in
\noindent{\small {\small\bf Lemma 4.3.} If  $u_{-}>u_{+},$  then we have
$$ \lim_{\gamma\rightarrow1}u_{\ast\gamma}=\lim_{\gamma\rightarrow1}(v _{\ast\gamma}+\beta t)=\lim_{\gamma\rightarrow1}\sigma_{1}^{\gamma}(t)=\lim_{\gamma\rightarrow1}\sigma_{2}^{\gamma}(t)=u_{\delta}(t),\eqno{(4.11)}
$$
and
$$ \lim_{\gamma\rightarrow1}\int^{\sigma_{2}^{\gamma}(t)}_{\sigma_{1}^{\gamma}(t)}\rho_{\ast\gamma}d\xi=u_{\delta}(t)[\rho]-[\rho (v+\beta t)]
=\frac{1}{2}(\rho_{-}+\rho_{+})(u_{-}-u_{+}), \eqno{(4.12)}
$$
where $u_{\delta}(t)=\frac{1}{2}(u_{-}+u_{+})+\beta t.$

\vskip 0.1in
\noindent{\small {\small\bf Proof.} It
follows from (1.2), (4.7), (4.8) and Lemma 4.2 that
$$\lim_{\gamma\rightarrow1}u_{\ast\gamma}=\lim_{\gamma\rightarrow1}(v_{\ast\gamma}+\beta t)=u_{-}+\beta t-\lim_{\gamma\rightarrow1}\sqrt{\frac{\frac{\gamma-1}{2}(\rho_{\ast\gamma}^{\gamma-1}-\rho_{-}^{\gamma-1})}
{(\rho_{\ast\gamma}+\rho_{-})(\rho_{\ast\gamma}-\rho_{-})}}(\rho_{\ast\gamma}-\rho_{-})
$$$$
=u_{-}+\beta t-\sqrt{a}=u_{-}+\beta t-\frac{1}{2}
(u_{-}-u_{+})=u_{\delta}(t),
$$
$$\lim_{\gamma\rightarrow1}\sigma_{1}^{\gamma}(t)=
u_{-}
+\beta t-\lim_{\gamma\rightarrow1}\rho_{\ast\gamma}\sqrt{\frac{\frac{\gamma-1}{2}(\rho_{\ast\gamma}^{\gamma-1}-\rho_{-}^{\gamma-1})}{(\rho_{\ast\gamma}
+\rho_{-})(\rho_{\ast\gamma}-\rho_{-})}}
=u_{-}+\beta t-\sqrt{a}=u_{\delta}(t),\eqno{}$$$$ \lim_{\gamma\rightarrow1}\sigma_{2}^{\gamma} (t)=\lim_{\gamma\rightarrow1}\bigg(v_{\ast\gamma}
+\beta t+\rho_{+}\sqrt{\frac{\frac{\gamma-1}{2}(\rho_{+}^{\gamma-1}-\rho_{\ast\gamma}^{\gamma-1})}{(\rho_{+}+\rho_{\ast\gamma})(\rho_{+}-\rho_{\ast\gamma})}}\,\,\bigg)
=\lim_{\gamma\rightarrow1}(v_{\ast\gamma}+\beta t)=u_{\delta}(t),\eqno{}$$
which immediately lead to $  \lim\limits_{\gamma\rightarrow1}u_{\ast\gamma}= \lim\limits_{\gamma\rightarrow1}\sigma_{1}^{\gamma}(t)=\lim\limits_{\gamma\rightarrow1}\sigma_{2}^{\gamma}(t)=u_{\delta}(t).
$

 From the first equations of the Rankine-Hugoniot conditions (3.8) for $S_{1}^{\gamma}$ and $S_{2}^{\gamma}$,   we have
 $$ \sigma_{1}^{\gamma}(t)(\rho_{-}-\rho_{\ast\gamma})=\rho_{-}(u_{-}+\beta t)-\rho_{\ast\gamma}(v_{\ast\gamma}+\beta t),\eqno{(4.13)}$$and $$\sigma_{2}^{\gamma}(t)(\rho_{\ast\gamma}-\rho_{+})=\rho_{\ast}(v_{\ast\gamma}+\beta t)-\rho_{+}(u_{+}+\beta t). \eqno{(4.14)}$$
From  (4.13),
(4.14) and (4.11), we get
 $$ \lim_{\gamma\rightarrow1}\rho_{\ast\gamma}(\sigma_{2}^{\gamma}(t)-\sigma_{1}^{\gamma}(t))=\lim_{\gamma\rightarrow1}(\rho_{-}(u_{-}+\beta t)
 -\sigma_{1}^{\gamma}(t)\rho_{-}+\sigma_{2}^{\gamma}(t)\rho_{+}-\rho_{+}(u_{+}+\beta t))
$$$$\,\,\,\,\,\,\,\,\,\,\,\,\,\,\,\,\,\,\,\,\,\,\,\,\,\,\,\,\,\,\,\,\,\,\,\,\,\,\,\,=u_{\delta}(t)[\rho]-[\rho (v+\beta t)]=\frac{1}{2}(\rho_{-}+\rho_{+})(u_{-}-u_{+}).\eqno{(4.15)}$$
Then, from (4.15), we obtain (4.12) immediately. The proof is completed. $~~\Box$

\vskip 0.1in
\noindent{\small {\small\bf Remark 4.1.} It can be concluded from Lemmas 4.2-4.3 that, when $\gamma\rightarrow1$,  the two shock  curves $S_{1}^{\gamma}$ and $S_{2}^{\gamma}$ will coincide,  the intermediate density $\rho_{\ast\gamma}$ becomes singular,  the limit of $\rho_{\ast\gamma}$ possesses a singularity which is a weighed Dirac delta function with the speed  $u_{\delta}(t)$.

\vskip 0.1in
\noindent{\small {\small\bf Remark 4.2.} It can be concluded from
 Lemma 4.3 that, when $\gamma\rightarrow1$, the velocities of two shocks $S_{1}^{\gamma}$ and $S_{2}^{\gamma}$ and the intermediate  $u_{\ast\gamma}$  of (1.1) approach to $u_{\delta}(t)$, which determines  the delta shock solution of  the pressureless Euler system with the
Coulomb-like friction term,  and the intermediate density $\rho_{\ast\gamma}$ between the two shocks tends to a weighted
$\delta$-measure which forms the delta shock.

\vskip 0.1in

From above analysis, we have the following result.
\vskip 0.1in

\vskip 0.1in
\noindent{\small {\small\bf Theorem 4.4.}
 For $u_{+}<u_{-}$, as $\gamma\rightarrow1$, the Riemann solution
containing two shocks  of (1.1) with the  Riemann initial data $(\rho_{\pm}, u_{\pm})$ constructed in Section
3 converges to
 a delta shock solution of system (1.4) with the same Riemann initial data $(\rho_{\pm}, u_{\pm})$.

\baselineskip 15pt
 \sec{\Large\bf 4.2.\quad    Formation of vacuum state for  system (1.1) }
In this subsection, we study the formation of vacuum state for the  Riemann solutions  containing two rarefaction waves of    (1.1) with the  Riemann initial data $(\rho_{\pm}, u_{\pm})$ as  $\gamma\rightarrow1$.
\vskip 0.1in
\noindent{\small {\small\bf Lemma 4.5.  }
  If  $u_{-}<u_{+}<u_{-}+2,$    then   there exists $\gamma_{1}>0 $  such  that
$(\rho_{+},u_{+})\in I(\rho_{-},u_{-})$   when  $1<\gamma<1+\gamma_{1}.$

\vskip 0.1in \noindent{\small {\small\bf Proof. }
It can be derived from (3.6) and (3.7)  that
  all possible states  $(\rho,v)$ that can be connected to  the left state $(\rho_{-},u_{-})$  on the right by a 1-rarefaction wave $R_{1}^{\gamma}$ or a 2-rarefaction wave $R_{2}^{\gamma}$ should satisfy
$$R_{1}^{\gamma}(\rho_{-},u_{-}):
 v+\rho^{\frac{\gamma-1}{2}}=u_{-}+\rho_{-}^{\frac{\gamma-1}{2}},\,\,\,\,v>u_{-},\,\,\,
    \rho<\rho_{-},
 \eqno{(4.16)}  $$
$$R_{2}^{\gamma}(\rho_{-},u_{-}):
 v-\rho^{\frac{\gamma-1}{2}}=u_{-}-\rho_{-}^{\frac{\gamma-1}{2}},\,\,\,\,\,\,v>u_{-},\,\,\,
    \rho>\rho_{-}.
\eqno{(4.17)}  $$

Similarly, it can be derived from  (3.7)  that
  all possible states  $(\rho,v)$ that can be connected to  the left state $(0,\widetilde{v}_{\ast}^{\gamma})$  on the right by  a 2-rarefaction wave $R_{2}^{\gamma}$ should satisfy
$$R_{2}^{\gamma}(0,\widetilde{v}_{\ast}^{\gamma}):
 v-\rho^{\frac{\gamma-1}{2}}=u_{-}+\rho_{-}^{\frac{\gamma-1}{2}},\,\,\,\,\,\,v>u_{-}+\rho_{-}^{\frac{\gamma-1}{2}},\,\,\,
    \rho>0.
\eqno{(4.18)}  $$

If $u_{-}<u_{+}<u_{-}+2$, $\rho_{+}\neq\rho_{-}$ and  $(\rho_{+}, u_{+})\in I(\rho_{-}, u_{-}), $     then we can see intuitively from  Figure 1 together with (4.16)-(4.18) that
$$
 u_{+}>u_{-}+\rho_{-}^{\frac{\gamma-1}{2}}-\rho_{+}^{\frac{\gamma-1}{2}},\,\,
    \rho_{+}<\rho_{-},
 \eqno{(4.19)}  $$
$$u_{+}>u_{-}-\rho_{-}^{\frac{\gamma-1}{2}}+\rho_{+}^{\frac{\gamma-1}{2}},\,\,
    \rho_{+}>\rho_{-},
\eqno{(4.20)}  $$
and $$u_{+}<u_{-}+\rho_{-}^{\frac{\gamma-1}{2}}+\rho_{+}^{\frac{\gamma-1}{2}},\,\,
    \rho_{+}>0.
\eqno{(4.21)}  $$ According to (4.19)-(21), we obtain that
$$|\rho_{-}^{\frac{\gamma-1}{2}}-\rho_{+}^{\frac{\gamma-1}{2}}|<u_{+}-u_{-}<\rho_{-}^{\frac{\gamma-1}{2}}+\rho_{+}^{\frac{\gamma-1}{2}}, \,\,\,\,\,\,\,\,\,\,
    \rho_{+}>0,\,\rho_{-}>0.$$
From
$  \lim\limits_{{\gamma\rightarrow1}}(\rho_{-}^{\gamma-1}-\rho_{+}^{\gamma-1})=0
  $ and $  \lim\limits_{{\gamma\rightarrow1}}(\rho_{-}^{\gamma-1}+\rho_{+}^{\gamma-1})=2>u_{+}-u_{-},
  $
it follows that there exists $\gamma_{1}>0$ small enough such that, when $ 1<\gamma <1 +\gamma_{1}$,  we have
$$|\rho_{-}^{\frac{\gamma-1}{2}}-\rho_{+}^{\frac{\gamma-1}{2}}|<u_{+}-u_{-}<\rho_{-}^{\frac{\gamma-1}{2}}+\rho_{+}^{\frac{\gamma-1}{2}}, \,\,\,\,\,\,\,\,\,\,
    \rho_{+}>0,\,\rho_{-}>0. \eqno{(4.22)}$$
  Then,  it is obvious
that $(\rho_{+}, u_{+})\in I( \rho_{-}, u_{-})$ when $ 1<\gamma <1 +\gamma_{1}$.
 The proof is completed. $~\Box$

\vskip 0.1in

When  $u_{-}<u_{+}<u_{-}+2,$  by Lemma 4.5,  for any given $\gamma\in (1, 1+\gamma_{1}),$  the Riemann solution  of (1.1)
with the  Riemann initial data $(\rho_{\pm}, u_{\pm})$ is as follows$$(\rho_{-}, u_{-}+\beta t)+R_{1}^{\gamma}+(\rho_{\ast\gamma}, v_{\ast\gamma}+\beta t)+R_{2}^{\gamma}+(\rho_{+}, u_{+}+\beta t),\eqno{(4.23)}
$$
where
$$R_{1}^{\gamma}:\,\,\left\{\begin{array}{ll} \frac{dx}{dt}=\lambda_{1}^{\gamma}=v+\beta t-\frac{\gamma-1}{2}\rho^{\frac{\gamma-1}{2}},
\\v+\rho^{\frac{\gamma-1}{2}}=u_{-}+\rho_{-}^{\frac{\gamma-1}{2}},
\,\,\,\,\rho_{\ast\gamma}\leq\rho\leq\rho_{-}, \end{array} \right.\eqno{(4.24)}$$ and
$$R_{2}^{\gamma}:\,\,\left\{\begin{array}{ll}  \frac{dx}{dt}=\lambda_{2}^{\gamma}=v+\beta t+\frac{\gamma-1}{2}\rho^{\frac{\gamma-1}{2}},
\\v-\rho^{\frac{\gamma-1}{2}}=u_{+}-\rho_{+}^{\frac{\gamma-1}{2}},
\,\,\,\,\rho_{\ast\gamma}\leq\rho\leq\rho_{+}.\end{array} \right.\eqno{(4.25)}$$
Thus,  from (4.24)  and (4.25), we can  derive  that
 $$ u_{+}-u_{-}=\rho_{+}^{\frac{\gamma-1}{2}}+\rho_{-}^{\frac{\gamma-1}{2}}-2\rho_{\ast\gamma}^{\frac{\gamma-1}{2}},\,\,\,\,
    \,\,\,\,\,\,\,\,\, \rho_{\ast\gamma}\leq\rho_{\pm}.
  \eqno{(4.26)}.$$  which  implies  the phenomenon of vacuum occurs as  $\gamma\rightarrow1$.

  \vskip 0.1in

  \vskip 0.1in\vskip 0.1in
\noindent{\small {\small\bf Theorem 4.6.}
Let $u_{-}<u_{+}< u_{-}+2$.   For any fixed $\gamma \in(1,2)$, assume that $(\rho_{\gamma}, u_{\gamma})(t,x)$ is a Riemann solution containing two rarefaction waves  of (1.1) with the  Riemann initial data $(\rho_{\pm}, u_{\pm})$ constructed in Section
3.  Then,  as $\gamma\rightarrow1$, the vacuum state occurs, and two
rarefaction waves become two contact discontinuities connecting the  states $(\rho_{\pm}, u_{\pm}+\beta t)$
and the vacuum $(\rho =0)$, which form a vacuum solution of system (1.4)  with the same initial data $(\rho_{\pm}, u_{\pm})$.

 \vskip 0.1in

\noindent{\small {\small\bf Proof.}  If $\lim\limits_{\gamma\rightarrow 1}\rho_{\ast\gamma}=K\in(0,\min(\rho_{-},\rho_{+})),$  then
taking the limit $\gamma\rightarrow 1$ in (4.27),  we have $u_{+}=u_{-}$, which contradicts with $u_{-}<u_{+}$. Thus $\lim\limits_{\gamma\rightarrow 1}\rho_{\ast\gamma}=0$, which means the vacuum occurs
as $\gamma\rightarrow 1$. Moreover, as $\gamma\rightarrow1$, one can directly derive from (4.24) and (4.25) that
$$\lim\limits_{\gamma\rightarrow 1}v=u_{-}  \,\,\mathrm{on}\,\, R_{1}^{\gamma},\,\,\,\,\,\,\,\,\,\,\,\,\,\,\,\lim\limits_{\gamma\rightarrow 1}v=u_{+} \,\,\mathrm{on}\,\, R_{2}^{\gamma},\eqno{(4.27)}  $$and
$$
\left\{
             \begin{array}{ll}
               \lambda_{1}^{\gamma}=\frac{\gamma+1}{2}v-\frac{\gamma-1}{2}u_{-}-\frac{\gamma-1}{2}\rho_{-}^{\frac{\gamma-1}{2}}
               +\beta t, \\
               \lambda_{2}^{\gamma}=\frac{\gamma+1}{2}v-\frac{\gamma-1}{2}u_{+}+\frac{\gamma-1}{2}\rho_{+}^{\frac{\gamma-1}{2}}
               +\beta t.
  \end{array}
\right.
\eqno{(4.28)}  $$
(4.27) and (4.28) imply that $$\lim_{\gamma\rightarrow 1}\lambda_{1}^{\gamma}=u_{-}+\beta t,  \,\,\,\,\,\, \,\,\,\,\,\lim_{\gamma\rightarrow 1}\lambda_{2}^{\gamma}=u_{+}+\beta t.\eqno{(4.29)}  $$
 The proof is completed. $~\Box$
\vskip 0.1in

\baselineskip 15pt
 \sec{\Large\bf 5.\quad  Numerical results for (1.1)}

 In this section, in order to  verify the validity of the formation of $\delta$-shocks  and vacuum states for system (1.1) mentioned in section 4, we present  two selected groups of representative numerical simulations. A number of iterative numerical trials are executed to guarantee what we demonstrate are not numerical objects. To discretize the system, we use the
  fifth-order weighted essentially non-oscillatory scheme and third-order Runge-Kutta method  [15, 29]  with the mesh 200 cells. The numerical simulations are consistent with the theoretical analysis.

\baselineskip 15pt
 \sec{\Large\bf 5.1.\quad    Formation of delta shock wave  }
 When $u_{+}<u_{-}$, we compute the solution of the Riemann problem of (1.1) with $\beta=2$  and
take the initial data as follows:
$$ (\rho, u)(0, x) =\left\{\begin{array}{ll} (1.5,\,
2),\,\,\,\,\,x< 0,\\(2,\,
-1),\,\,\,\,x> 0.\end{array} \right.\eqno{(5.1)}$$
The numerical simulations for different choices of  $\gamma$
 ( $\gamma=1.7$, $1.05$, $1.001$,
 and the time $t=0.2$ )
are presented in Figs. 2-4 which show the process of concentration  and formation of the delta
shock wave in the pressureless limit of solutions containing two shocks.

\begin{figure}[htbp]
\centering
\begin{minipage}[c]{0.45\textwidth}
\centering
\includegraphics[width=2.5in]{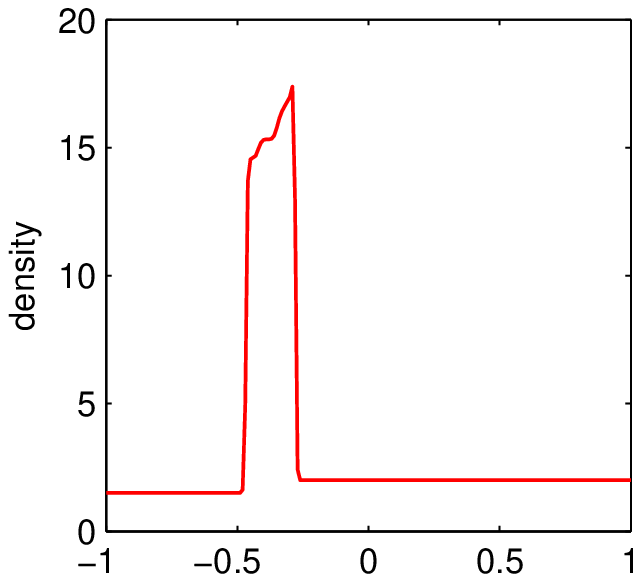}
\end{minipage}%
\begin{minipage}[c]{0.45\textwidth}
\centering
\includegraphics[width=2.5in]{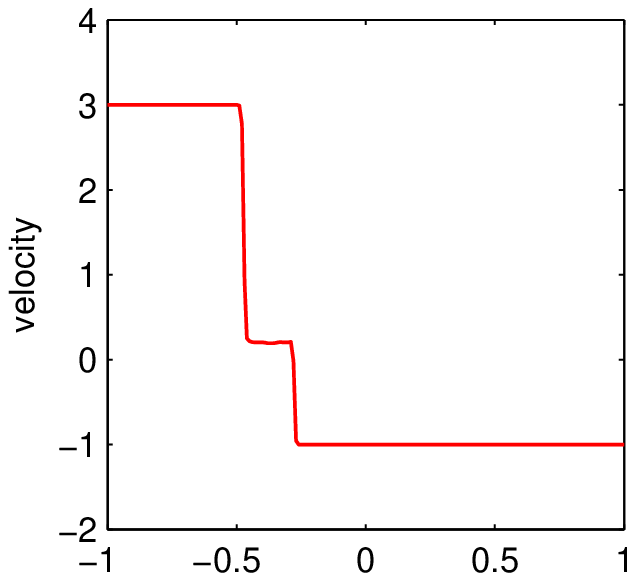}
\end{minipage}%
\end{figure}

\centerline{\bf Fig. 2.\,\,   Density (left) and velocity (right) for $\gamma=1.7$.}

\begin{figure}[htbp]
\centering
\begin{minipage}[c]{0.45\textwidth}
\centering
\includegraphics[width=2.5in]{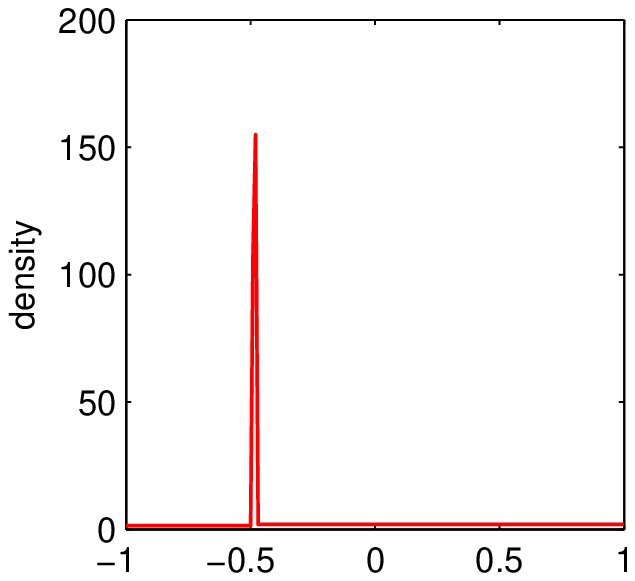}
\end{minipage}%
\begin{minipage}[c]{0.45\textwidth}
\centering
\includegraphics[width=2.5in]{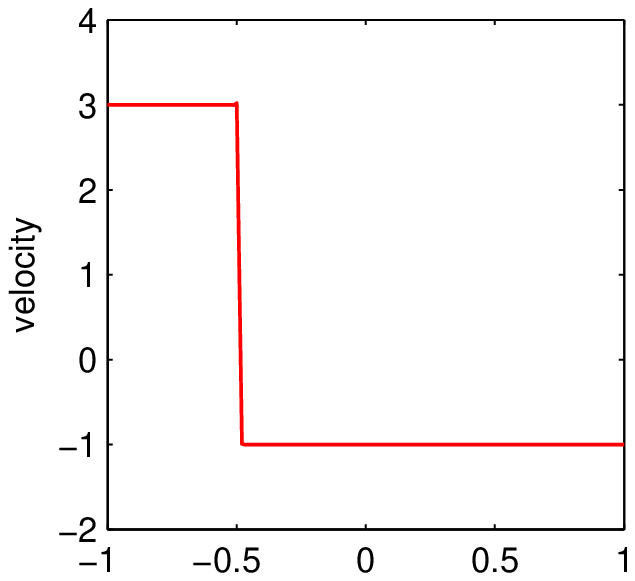}
\end{minipage}%
\end{figure}

\centerline{\bf Fig. 3.\,\,   Density (left) and velocity (right) for $\gamma=1.05$.}

\newpage
\begin{figure}[htbp]
\centering
\begin{minipage}[c]{0.45\textwidth}
\centering
\includegraphics[width=2.5in]{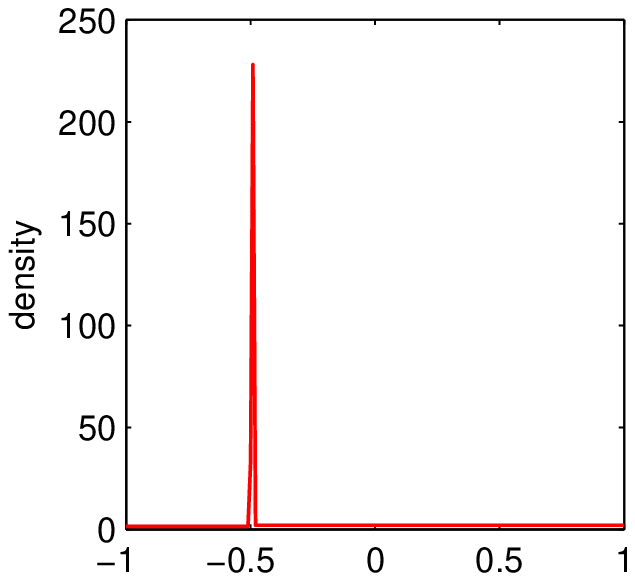}
\end{minipage}%
\begin{minipage}[c]{0.45\textwidth}
\centering
\includegraphics[width=2.5in]{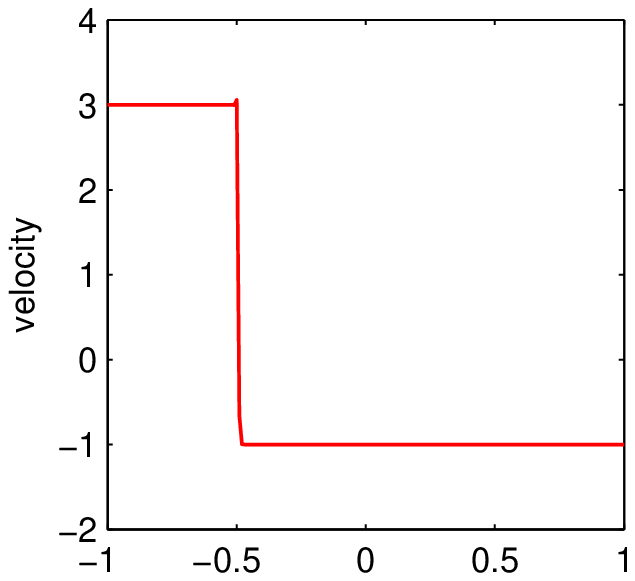}
\end{minipage}%
\end{figure}

\centerline{\bf Fig. 4.\,\,   Density (left) and velocity (right) for $\gamma=1.001$.}

We can clearly see from these numerical results that, as $\gamma$ decreases, the locations of the two shocks become closer and closer, and the density of the intermediate state increases dramatically, while  the velocity becomes a piecewise constant function. Finally, as $\gamma$ tends to one, along with the intermediate state, the two shocks coincide to form  the delta shock wave of  the pressureless Euler system with the
Coulomb-like friction term (1.4), while the velocity keeps a step function. The numerical simulations are in complete agreement with the theoretical analysis in section 4.1.

\baselineskip 15pt
 \sec{\Large\bf 5.2.\quad    Formation of the vacuum state  }
 When $u_{-}<u_{+}$, we compute the solution of the Riemann problem of (1.1) with $\beta=2$ and
take the initial data as follows:
$$ (\rho, u)(0, x) =\left\{\begin{array}{ll} (1,
-0.1),\,\,\,x< 0,\\(4,
1),\,\,\,\,\,\,\,\,\,\,\,x> 0.\end{array} \right.\eqno{(5.2)}$$
 The numerical simulations for different choices of $\gamma$
 ($\gamma=1.8$, $1.18$, $1.01$
and the time $t=0.2$), are
presented in Figs. 5-7 which show the process of cavitation and formation of the vacuum state in the
pressureless limit of solutions containing two rarefaction waves.

From these numerical results, we can clearly observe that,  when $\gamma$ decreases, the boundaries of two
rarefaction waves become closer and closer, along with the intermediate state, the density tends to zero,
while the velocity becomes a linear function. In the end, as $\gamma$ tends to one,  a two-rarefaction-wave solution tends to  a two-contact-discontinuity solution with a vacuum state of  the pressureless Euler system with the
Coulomb-like friction term  (1.4). The numerical simulations are in complete agreement with the theoretical analysis in section 4.2.

\newpage
\begin{figure}[htbp]
\centering
\begin{minipage}[c]{0.45\textwidth}
\centering
\includegraphics[width=2.5in]{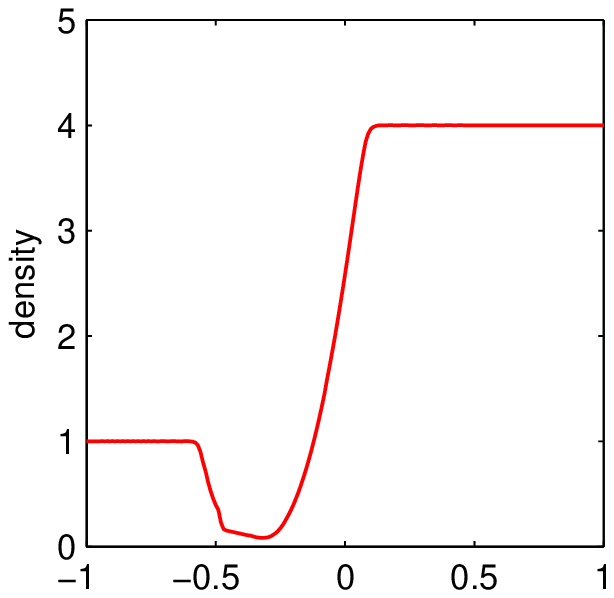}
\end{minipage}%
\begin{minipage}[c]{0.45\textwidth}
\centering
\includegraphics[width=2.5in]{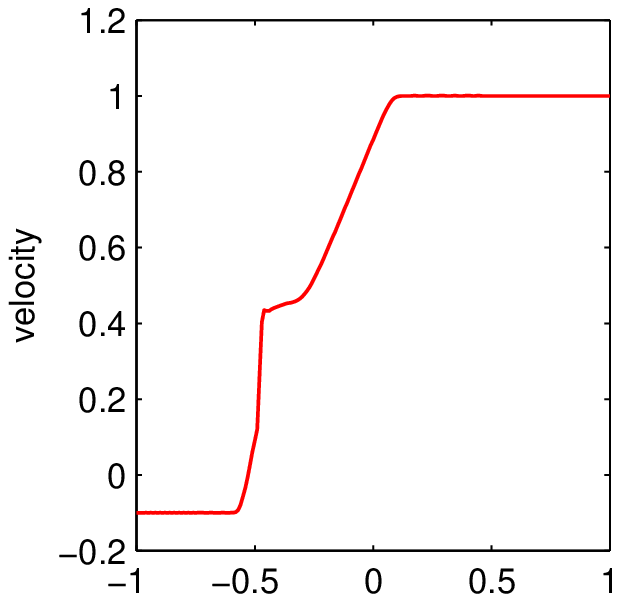}
\end{minipage}%
\end{figure}

\centerline{\bf Fig. 5.\,\,   Density (left) and velocity (right) for $\gamma=1.8$.}

\begin{figure}[htbp]
\centering
\begin{minipage}[c]{0.45\textwidth}
\centering
\includegraphics[width=2.5in]{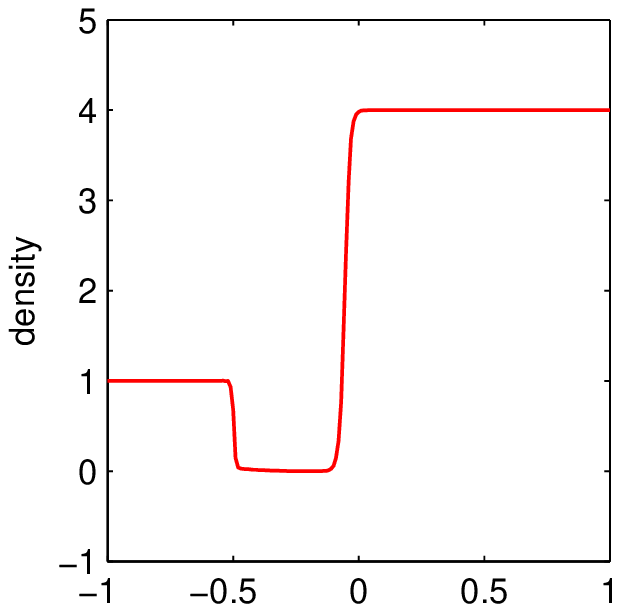}
\end{minipage}%
\begin{minipage}[c]{0.45\textwidth}
\centering
\includegraphics[width=2.5in]{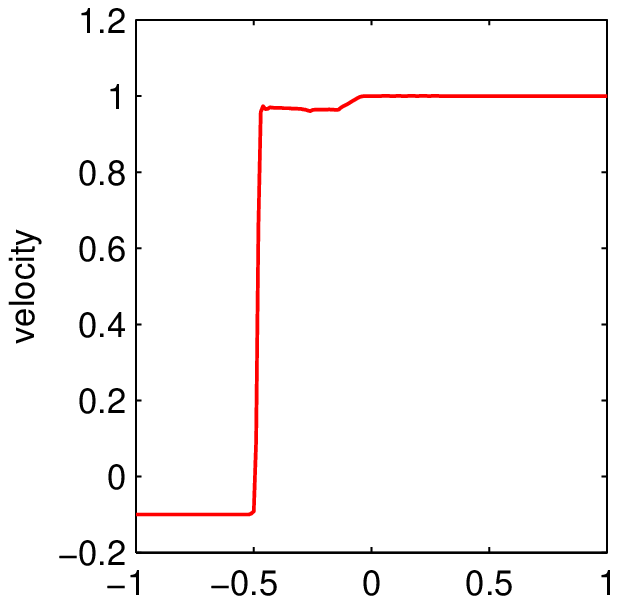}
\end{minipage}%
\end{figure}

\centerline{\bf Fig. 6.\,\,   Density (left) and velocity (right) for $\gamma=1.18$.}

\begin{figure}[htbp]
\centering
\begin{minipage}[c]{0.45\textwidth}
\centering
\includegraphics[width=2.5in]{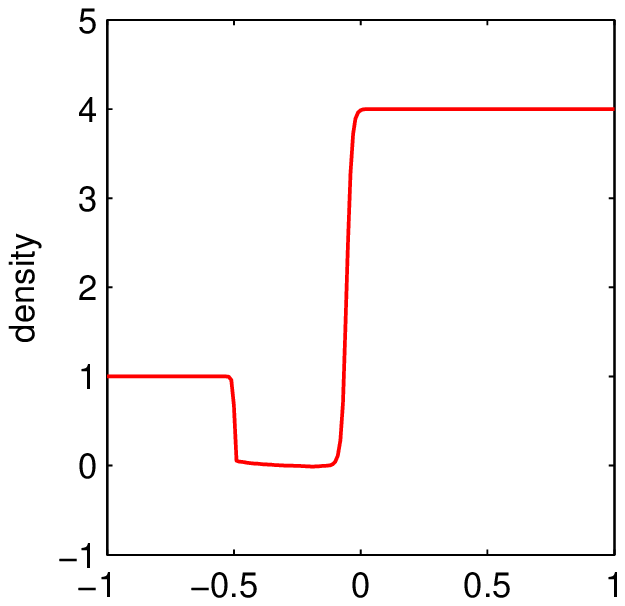}
\end{minipage}%
\begin{minipage}[c]{0.45\textwidth}
\centering
\includegraphics[width=2.5in]{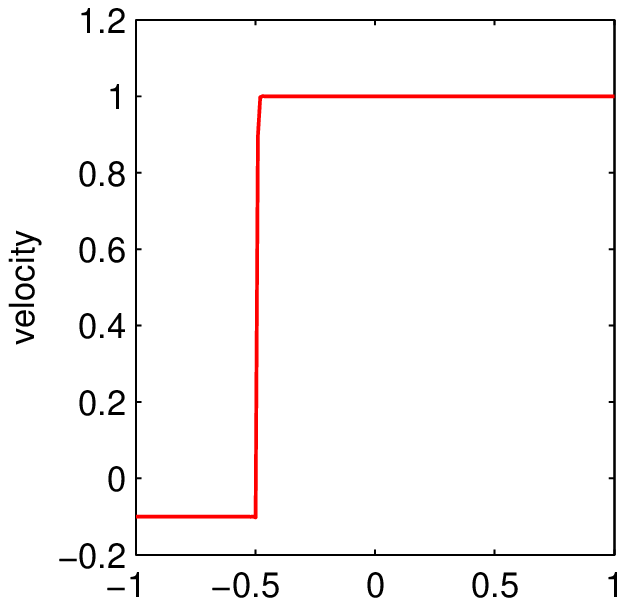}
\end{minipage}%
\end{figure}

\centerline{\bf Fig. 7.\,\,   Density (left) and velocity (right) for $\gamma=1.01$.}

\end{document}